\documentclass[final,5p,11pt,review]{elsarticle}
\usepackage{amssymb, amsmath}
\usepackage{amsthm}
\usepackage{bm}

\theoremstyle{definition}

\theoremstyle{remark}

\usepackage[section]{placeins} 
\usepackage[colorlinks]{hyperref}
\hypersetup{
     colorlinks   = true,
     citecolor    = blue
}
\usepackage{lineno}
\usepackage[section]{placeins} 
\usepackage{lineno}

\usepackage{setspace}
\usepackage{amssymb}
\usepackage{color}
\usepackage{colortbl} 
\usepackage{xcolor}
\usepackage{capt-of}
\usepackage{savesym}
\savesymbol{AND}
\usepackage{algorithmic}
\usepackage{algorithm}

\usepackage{setspace}
\journal{Engineering Analysis with Boundary Elements}

\begin{document}
\begin{frontmatter}

\title{An improved radial basis-pseudospectral method with hybrid Gaussian-cubic kernels}
\author[label1]{Pankaj K Mishra}
\address[label1]{Advanced Computational Seismology Laboratory, Indian Institute of Technology, Kharagpur, India}
\ead{pankajkmishra01@gmail.com}

\author[label1]{Sankar K Nath\corref{cor1}}
\cortext[cor1]{Corresponding Author}
\ead{nath@gg.iitkgp.ernet.in}
\author[label4]{Gregor Kosec}
\address[label4]{Parallel and Distributed Computing Laboratory, Jozef Stefan Institute, Slovenia}
\ead{gkosec@ijs.si}
\author[label5]{Mrinal K Sen}
\address[label5]{Institute for Geophysics, University of Texas at Austin, USA}
\ead{msentx@gmail.com}

\begin{abstract}
While pseudospectral (PS) methods can feature very high accuracy, they tend to be severely limited in terms of geometric flexibility. Application of global radial basis functions overcomes this, however at the expense of problematic conditioning (1) in their most accurate flat basis function regime, and (2) when problem sizes are scaled up to become of practical interest. The present study considers a strategy to improve on these two issues by means of using hybrid radial basis functions that combine cubic splines with Gaussian kernels. The parameters, controlling Gaussian and cubic kernels in the hybrid RBF, are selected using global particle swarm optimization. The proposed approach has been tested with radial basis-pseudospectral method for numerical approximation of Poisson, Helmholtz, and Transport equation. It was observed that the proposed approach significantly reduces the ill-conditioning problem in the RBF-PS method, at the same time, it preserves the stability and accuracy for very small shape parameters. The eigenvalue spectra of the coefficient matrices in the improved algorithm were found to be stable even at large degrees of freedom, which mimic those obtained in pseudospectral approach. Also, numerical experiments suggest that the hybrid kernel performs significantly better than both pure Gaussian and pure cubic kernels.
\end{abstract}

\begin{keyword}
radial basis function, pseudospectral method, ill-conditioning, partial differential equations 
\end{keyword}

\end{frontmatter}
\onecolumn 
\section{Introduction}
\noindent Pseudospectral (PS) methods are highly accurate and higher-order numerical methods, which use polynomials as basis functions. In two or higher dimensions, PS method tend to be limited in terms of geometric flexibility \citep{Forn1996}. A typical variant of PS methods is Chebyshev pseudospectral method (CHEB-PS), which uses Chebyshev polynomials as basis functions. In order to make the PS method  geometrically flexible, Fasshauer \citep{Fass06} proposed the application of infinitely smooth radial basis functions (RBFs) in pseudospectral formulation and interpreted the combined approach as meshless radial basis-pseudospectral (RBF-PS) method. Gaussian RBF is one such infinitely smooth RBF, which results in a positive definite system ensuring uniqueness in the interpolation. It is often found in the application of smooth RBFs that scaling the radial kernel by reducing the shape parameter to a smaller value, \textit{i.e.}, making it ``flat'' reduces the error in the approximation, as the ``flat'' limit of infinitely smooth RBF converges to a polynomial interpolant \citep{Driscoll2002413,Fornberg200437}. Larsson and Fornberg \citep{Larsson2005103} have shown that it is possible to get even more accurate results with Gaussian RBF in the ``flat'' range, \textit{i.e}, just before it converges to polynomial interpolants. Although global RBF methods are relatively costly because of the full and dense matrices arising in the linear system, their accuracy and convergence makes them desirable, especially for problems in solid mechanics. In recent years, RBF-PS method has been effectively applied to computational mechanics \citep{Ferreira2006134,Ferreira2007202,Krowiak2016}, nonlinear equations \citep{Uddin2013619}, and thermal convection in 3D spherical shells \citep{GGGE1704}, etc. Application of an infinitely smooth RBF in pseudospectral mode, however, brings an inherent limitation, as the global approximation of RBFs gets severely ill-conditioned at higher degrees of freedom as well as at low shape parameters. Such limitations constraint the well-posedness of the RBF-PS algorithm only to few nodes in the domain with relatively large shape parameter range. Typical quantification of such limitations can be found in \citep{Fass06}, where the RBF-PS algorithm was found to be well-posed upto $24\times24$ nodes for 2D Helmholtz's equation and $18$ nodes for 1D transport equation.

To deal with the ill-conditioning in RBF interpolation, Kansa and Hon {\color{blue}\citep{KansaHon2002}} performed numerical tests using various tools, \textit{viz.}, block partitioning or LU decomposition, matrix preconditioners, variable shape parameters, multizone methods, and node adaptivity. Other major contributions to deal with the mentioned problem are: a direct solution approach \citep{Cheng2005}, accelerated iterated approximate moving least squares {\color{blue}\citep{Fass2009}}, random variable shape parameters {\color{blue}\citep{Sarra20091239}}, Contour-Pad\'e  and RBF-QR algorithms  {\color{blue}\citep{Forn2011}}, series expansion {\color{blue}\citep{Fass2012}}, and regularized symmetric positive definite matrix factorization \citep{Sarra2014}, RBF-GA {\color{blue}\citep{Fornberg200760}}, Hilbert-Schmidt SVD {\color{blue}\citep{Fass2015}}, Weighted SVD {\color{blue}\citep{DeMarchi20131}}, use of Laurent series of the inverse of the RBF interpolation matrix \citep{Kindelan2016}, and RBF-RA \citep{Wright2017}, \textit{etc}. An alternative approach is radial basis finite difference (RBF-FD) method, which is a local version of RBF-PS method \citep{Chandhini2007,Bayona2010,Flyer2012,Flyer2016}. The only significant difference between RBF-PS and RBF-FD implementation is that instead of using all the nodes, later uses only few neighbour nodes for construction of differential matrices.

Recently, Mishra et. al. \cite{1512.07584} proposed novel radial basis functions by hybridizing Gaussian and cubic kernels, which could significantly reduce the ill-conditioning problem in scattered data interpolation. This hybrid kernel utilizes optimal proportion of the Gaussian and cubic kernel, which correspond to the defined optimization criterion. In this paper, we propose a well-conditioned radial basis-pseudospectral scheme for numerical approximation of PDEs, by incorporating hybrid Gaussian-cubic kernels as basis functions. We establish both the convergence and stability of this improved scheme, through several numerical examples including numerical approximation of time-independent and time-dependent PDEs. Hereafter, in this work, we will call this improved approach as hybrid radial basis function-pseudospectral approach (HRBF-PS).

Rest of the paper is structured as follows. We introduce the hybrid RBF in section \ref{sec:HGCR}, and the global particle swarm optimization algorithm for selecting the parameters of this hybrid RBF in section \ref{sec:pso}. Construction of differentiation matrices, and the RBF-PS scheme for numerical solution of PDEs have been explained in section \ref{sec:rbfps}. Finally we perform numerical tests by solving Poisson, Helmholtz, and transport equations using the improved RBF-PS method and exhibit the improvements, observed due to hybrid RBF over Gaussian and Cubic RBFs, in section \ref{sec:test}, followed by the conclusion. In appendix A, we explain the particle swarm optimization algorithm and its application in the contexts of numerical solution of PDEs with meshless methods.

\section{Hybrid Gaussian-cubic RBF}
\label{sec:HGCR}
\noindent Radial basis functions were proposed by Hardy \citep{JGR:JGR12292} for fitting topography on irregular surfaces using linear combination of a single symmetric basis functions, which was later found to have better convergence than many available approaches for interpolation \citep{Franke1979}. Some commonly used RBFs have been listed in Table (\ref{tab:rbf}). First application of RBFs for numerical solution of differential equations was proposed by Edward Kansa in 1990 \citep{Kansa1990127}. Since RBFs do not require to be interpolated on regular tensor grids, Kansa's method did not require``mesh'', therefore, it was termed as a meshless method. Infinitely smooth RBFs like Gaussian have been proven to provide invertible system matrix in such meshless methods. However, for small shape parameters, as well as large number of nodes in the domain, Gaussian RBF leads to solving an ill-conditioned system of equations. Cubic RBFs on the other hand, are finitely smooth radial basis functions, which, unlike Gaussian RBF, do not have any shape parameter. However, use of cubic RBF for shape function interpolation in meshless methods involves the risk of getting a singular system, for certain node arrangements. Recently, a hybrid RBF \citep{1512.07584}, by combining Gaussian and the cubic kernels, has been proposed which could utilize certain features of both the RBFs depending on the problem type under consideration, as given by 

\begin{eqnarray}
\phi(r) = \alpha e^{ -(\epsilon r)^2} + \beta r^3,
\end{eqnarray}

\noindent where, $\epsilon$ is the shape parameter of the radial basis function, which is a relatively new notation for the same. One advantage of using this new conventions is that all the RBFs depend on the shape parameter in a similar manner. It should be noted that there is another parallel convention for the shape parameter, which is commonly represented as `$c$' \citep{Chen2014}.  The conversion from old to new convention can be done by setting $c^2=1/\epsilon^2$ \citep{Fassbook2007}. The weight coefficients $\alpha$ and $\beta$ control the contribution of Gaussian and cubic kernel in the proposed hybridization depending upon the problem type.
\begin{table}[!htbp]
  \centering \footnotesize
  \begin{tabular}{lc}
    \hline
    Kernel & Mathematical expression \\
    \hline
    Multiquadric (MQ)                       & $ (1+(\epsilon r)^2)^{1/2}$ \\
    Inverse multiquadric (IMQ)              & $ (1+(\epsilon r)^2)^{-1/2} $\\
    Gaussian (GA)                             & $ e^{-(\epsilon r)^2}$\\
    Polyharmonic Spline (PHS)                & $\begin{cases} r^m ln(r) \qquad m =2,4,6,... \\ r^m \qquad\qquad m= 1, 3, 5,...\end{cases}
$ \\
    Wendland's (Compact Support)            & $(1-\epsilon r)^{4}_{+}(4\epsilon r+1)$\\
    \hline
\end{tabular}
 \caption{Some frequently used radial basis functions (radial kernels) and their mathematical expressions.}
 \label{tab:rbf}
\end{table}

\section{Parameter Optimization}
\label{sec:pso}
\noindent Since the shape parameter affects both the accuracy and stability of algorithms involving RBFs, finding its optimal value has been a critical issue in radial basis interpolation and its application in meshless methods \citep{Roque2010,Huang2010,Cheng2012}. The hybrid kernel, presented in this study, contains three parameters, \textit{i.e.}, $\epsilon$, $\alpha$, and $\beta$, an optimal combination of which will ensure the optimum convergence and stability of the associated algorithm. Particle swarm optimization (PSO) is a frequently used algorithm to decide the shape parameter in RBF network and its application in machine learning algorithm \citep{Liu2010,Esma2009}, however in context of numerical approximation of PDEs with meshless methods, it is generally decided with ad-hoc methods like solving the problem with various values of the shape parameter and visualizing the root mean square (RMS) error against it. This approach works only if the exact solution of the problem is known, which in practical cases, is often unknown. For such cases, in the context of scattered data interpolation, Rippa \citep{Rippa} proposed a statistical approach using leave-one-out-crossvalidation (LOOCV), which later got generalized for numerical solution of PDEs, by Fasshauer \citep{FASS20077}. Here we use a global particle swarm optimization algorithm, to decide the optimal values of the parameters of the hybrid kernel. We test two different objective functions: (1) RMS error, when the exact solution is known and (2) LOOCV criterion, when the exact solution is not known. Algorithm (\ref{alg:test}), explains the process of computing the objective function using LOOCV. Here $c_k$ is the $k^{th}$ coefficient for the interpolant on ``full data" set and $\mathbf{A}^{-1}_{kk}$ is the $k^{th}$ diagonal element in the inverse of the interpolation matrix for ``full data". A detailed discussion about the application of particle swarm optimization in this context has been given in Appendix A.

\begin{algorithm}[!htbp]
\begin{algorithmic}[1]
\STATE Fix a set of parameters $[\epsilon, \alpha, \beta]$
\FOR{all the N collocation points, \textit{i.e.},$k=1,...,N$}
\IF{using Rippa's simplified approach \citep{Rippa}}
 \STATE Compute the error vector $e_k$ as
\begin{eqnarray}
e_k =  \frac{c_k}{\mathbf{A}^{-1}_{kk}}.
\end{eqnarray}
  \ELSE
\STATE Compute the interpolant by excluding the $k^{th}$ point as following, (see equation (\ref{approximation}))
\begin{eqnarray}
\mathcal{I}(\bm{x}) = \sum_{j=1}^{N-1} c^{[k]}_{j} \phi (\parallel \bm{x}-\bm{x}^{[k]}_{j}\parallel).
\end{eqnarray}
\STATE Compute the $k^{th}$ element of the error vector $e_k$ 
\begin{eqnarray}
e_k = \mid\mathcal{I}(\bm{x}_k) - \mathcal{I}^{[k]}(\bm{x}_k) \mid,
\end{eqnarray}
  \ENDIF
\ENDFOR
\STATE Assemble the ``cost vector'' as $\bm{e} = [ e_1,..., e_N]^T$.
\STATE \noindent The optimization problem here, can be written in the mathematical form as, 
\[ Minimize \rightarrow \parallel\bm{e}\parallel (\epsilon, \alpha, \beta),  \]
subject to the following constrains, 
\[\epsilon \geq 0,\]
\[0 \leq \alpha \leq 1,\] 
\[0 \leq \beta \leq 1.\]
\end{algorithmic}
\caption{LOOCV for computing the objective function for parameter optimization. This algorithm uses the interpolation matrix, which is computed to construct various differentiation matrices in RBF-PS.}
\label{alg:test}
\end{algorithm} 

\section{RBF-PS Scheme}
\label{sec:rbfps}
\noindent Kansa's collocation method \citep{Kansa19901} is a frequently used approach for numerical solution of PDEs via meshless approach, which provides a solution as a continuous function. Computing such a continuous solution may make the algorithm relatively expensive, especially for time dependent problems \citep{Ferreira2006134}. Pseudospectral methods on the other hand, provide the solution on certain specified nodes. PS methods are often implemented by constructing differentiation matrices and substituting them into the involved differential equation. We briefly discuss the construction of differentiation matrices for a univariate case below, however, detailed explanations and implementation can be found in Nick Trefethen's book \citep{Trefethen2000}.

\noindent Let us assume a computational domain $\Omega$ and discretize it using $N$ number of nodes $\bm{x}_k, k=1,...,N$. The approximation $\bm{u}$ of the unknown field in a typical PDE can be written as a linear combination of some unknown coefficients $c_k$ and RBFs $\phi_k$ as given by
\begin{eqnarray}
\label{approximation}
u(\bm{x}_i)=\sum_{k=1}^{N} c_k \phi \left(\parallel \bm{x}_i-\bm{x}_k \parallel\right).
\end{eqnarray}
\noindent where $\phi_k = \phi(\parallel x - x_k \parallel)$. Equation (\ref{approximation}) can be written in the matrix form as 
\begin{eqnarray}
\label{approxmatrix}
\bm{u} = \bm{A}\bm{c},
\end{eqnarray}
where $A_{ik}= \phi \left( \parallel \bm{x}_i- \bm{x}_k \parallel \right)$ are the basis functions at nodes, and $\bm{c}=[c_1,...,c_N]^T$ are the corresponding unknown coefficients. The derivative of $\bm{u}$ can be computed by differentiating the basis functions,\textit{i.e.},
\begin{eqnarray}
\frac{d}{d\bm{x}_i}u(\bm{x}_i) = \sum_{k=1}^{N} c_k \frac{d}{d\bm{x}_i} \phi \left( \parallel \bm{x}_i- \bm{x}_k \parallel \right).
\end{eqnarray}
The matrix-vector notation of the above equation for the derivatives at collocation points $\bm{x}_i$ can be written as 
\begin{eqnarray}
\label{dm}
\bm{u}' = \bm{A}_x\bm{c},
\end{eqnarray}
where the derivative matrix $\bm{A}_x$ has the entries of derivatives of the radial basis functions, \textit{i.e.}, $\frac{d}{d\bm{x}_i}\phi \left( \parallel \bm{x}_i- \bm{x}_k \parallel \right)$. In the context of RBF-PS methods, the matrix $\bm{A}$ is radial basis interpolation matrix. Since Gaussian RBF is positive definite, its application in RBF-PS method ensures the invertibility of the interpolation matrix for distinct nodes. Moreover, the flat limit of Gaussian RBF converges to polynomial interpolants which can achieve spectral accuracy. It is well known that the global approximation with Gaussian RBF leads to solving severely ill-conditioned system, therefore, the well-posedness of the interpolation matrix in RBF-PS algorithm is a critical issue. Mishra et. al. \citep{1512.07584} have shown the well-posedness as well as the accuracy of hybrid Gaussian cubic kernels for the interpolation problem in the ``flat" region. In order to show the similar improvements in RBF-PS method, we will use this hybrid RBF in all the numerical tests.

\noindent Substituting the coefficient vector $\bm{c}$ from equation (\ref{approxmatrix}) to equation (\ref{dm}), we get
\begin{eqnarray}
\label{dm2}
\bm{u}' = \bm{A}_x \bm{A}^{-1}\bm{u}.
\end{eqnarray}
Hence, the corresponding differentiation matrix can be written as 
\begin{eqnarray}
\label{dm3}
\bm{D} = \bm{A}_x \bm{A}^{-1}.
\end{eqnarray}
Higher order derivatives and other complex linear operators can be computed by following the similar procedure, unlike the PS method, where higher order differentiation matrices can be computed as products of first order differentiation matrix \citep{Fass06}. For example, a typical linear operator $\mathcal{L}$ can be constructed as,
\begin{eqnarray}
\label{opL}
\mathcal{L} = \bm{A}_{\mathcal{L}} \bm{A}^{-1},
\end{eqnarray}
where element of the matrix $\bm{A}_{\mathcal{L}}$ are $\left[A_{\mathcal{L}}\right]_{ik}=\mathcal{L}\phi \left( \parallel \bm{x}_i- \bm{x}_k \parallel \right)$.

\noindent The boundary conditions can be incorporated in a RBF-PS scheme via two different approaches. In some cases, like periodic problems, the chosen basis function satisfies the boundary conditions. In other cases, the boundary conditions are explicitly enforced into the system. We explain this explicit enforcing of the boundary condition by considering the following linear elliptic differential equation
\begin{eqnarray}
\label{lineareq}
\mathcal{L}u(\bm{x}) = f(\bm{x}) \qquad \bm{x}\in\Omega,
\end{eqnarray}
with boundary conditions
\begin{eqnarray}
u(\bm{x}) = g(\bm{x}) \qquad \bm{x}\in \partial \Omega.
\end{eqnarray}
The discretized differential oprator without incorporating the boundary condition can be computed according to the equation (\ref{opL}). In order to incorporate the boundary conditions in the discretized operator, we replace the corresponding rows
of the $\mathcal{L}$ corresponding to the boundary collocation points by unit vectors
in the diagonal position and zeros at every other position, and then replace
the corresponding rows of $\bm{f}$ on the right-hand side by $\bm{g}$. Thus, the matrix-vector form of this problem can be is written as 
\begin{eqnarray}
\mathcal{L}_{bc}\bm{u} = 
\begin{bmatrix}
\bm{f} \\
\bm{g}
\end{bmatrix},
\end{eqnarray}
where $\mathcal{L}_{bc} = \begin{bmatrix}
\bm{A}_{\mathcal{L}} \\
\bm{A}
\end{bmatrix} 
\bm{A}^{-1}
$ is the modified discrete differential operator, which also contains the differential operators corresponding to the boundary conditions. A detailed discussion of this process can be found in \citep{Fass06}.
\section{Numerical Tests}
\label{sec:test}
\noindent In this section, we perform efficacy test of the improved RBF-PS method by considering several numerical examples adapted from Nick Trefethen \citep{Trefethen2000}.

\subsection{Poisson equation}
\label{s5.1}
\begin{figure}[hbtp]
\centering
\includegraphics[scale=0.4]{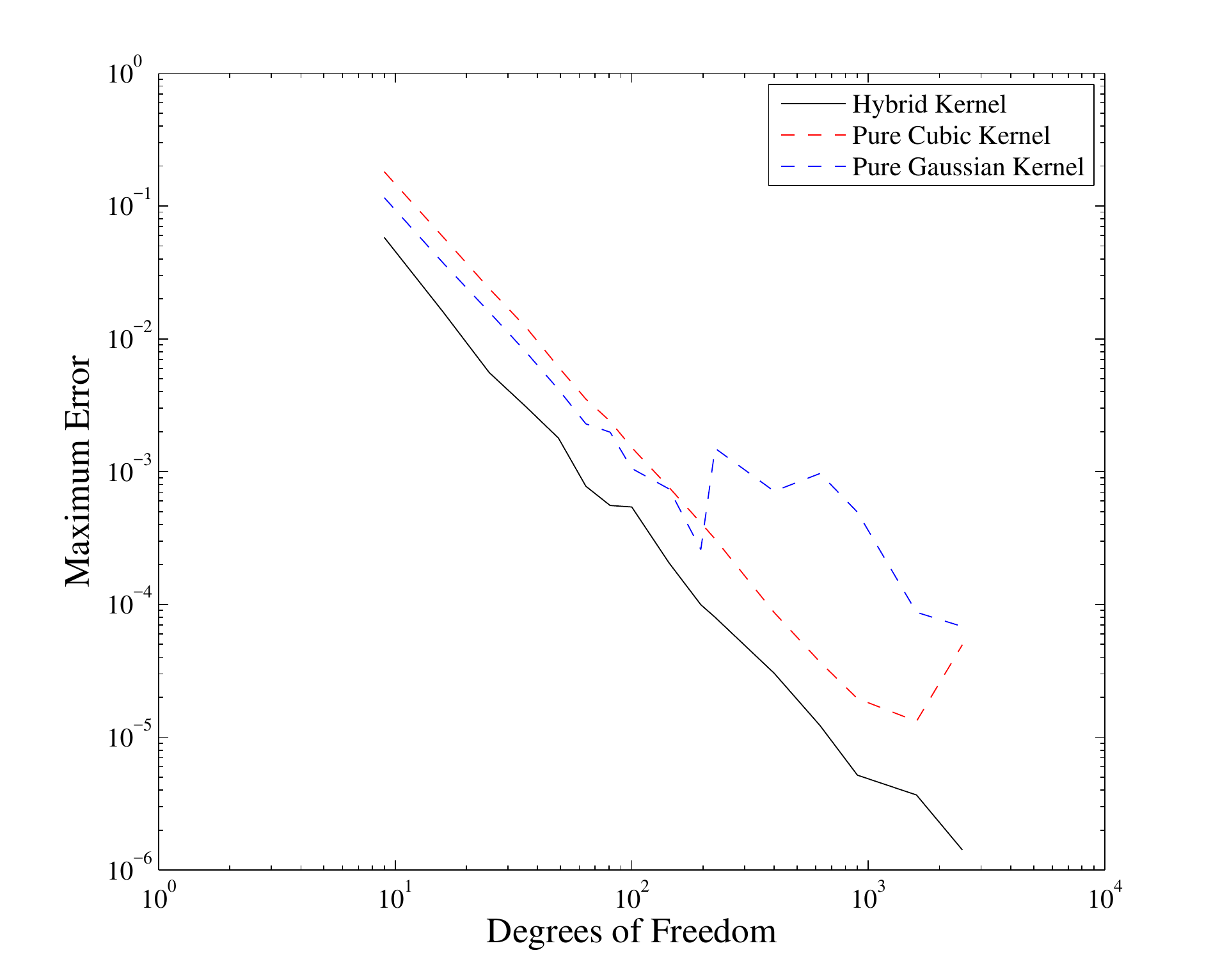}
\caption{The convergence pattern of the numerical approximation of Poisson equation with pure Gaussian, pure cubic, and hybrid kernel in RBF-PS scheme.}
\label{FNT1}
\end{figure}

\begin{figure}[hbtp]
\centering
\includegraphics[scale=0.4]{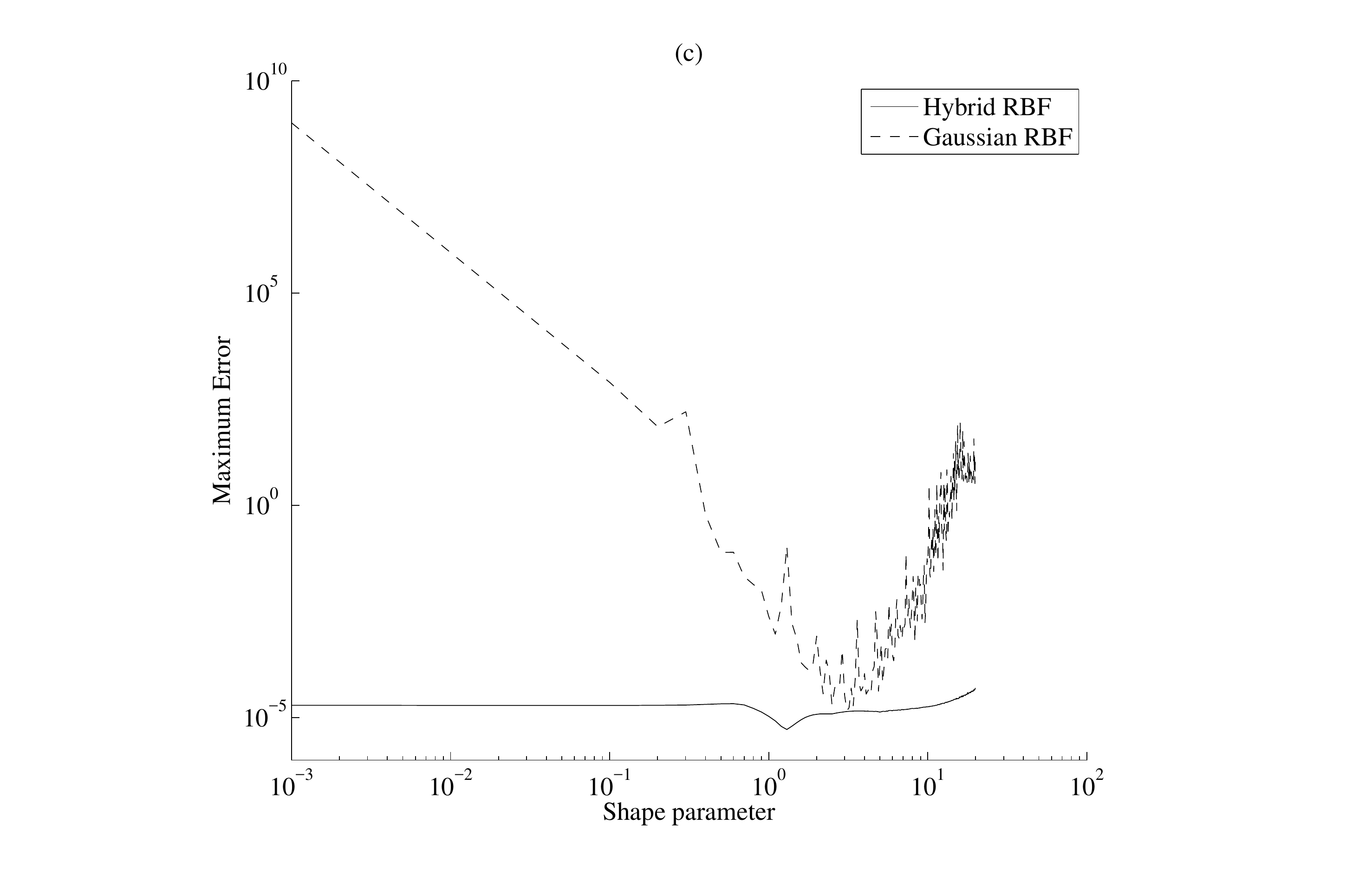}
\caption{Variation of the maximum error with different values of shape parameter for the Gaussian and the Hybrid kernel with optimized weight coefficients.}
\label{FNT122}
\end{figure}

\noindent To start with, we consider a simple linear and univariate boundary value problem in the domain $\Omega = [-1,1]$ with null Dirichlet boundary conditions, which is expressed as 
\begin{eqnarray}
 \frac{\partial^2 u}{\partial x^2} = e^{4x}, \qquad x\in \Omega.
 \label{possioneq}
\end{eqnarray}
The analytical solution of this problem is given as
\[u(x) = \left[ e^{4x}-tsinh(4)-cosh(4)\right].\]
We solve equation (\ref{possioneq}) with RBF-PS approach using hybrid Gaussian-cubic RBF, over Chebyshev grid points in the domain. Particle swarm optimization has been used to select the shape parameter $\epsilon$, and the weight coefficients $\alpha$, and $\beta$. Table \ref{NT1} contains the optimized values of the parameters and maximum error in this test for various degrees of freedom, \textit{viz.}, $[N=9,...,2500]$. Figure \ref{FNT1} exhibits excellent convergence of HRBF-PS algorithm used for this problem. Unlike Gaussian kernel, which becomes ill-conditioned, the hybrid kernel maintains the accuracy even with very small shape parameter, as shown in Figure \ref{FNT122}. Another implied advantage of using the hybrid kernel is that it can perform computation with relatively larger degrees of freedom without any special consideration, as it is well-posed.
  \begin{table*}[htbp]
     \centering \footnotesize
     \begin{tabular*}{\textwidth}{l@{\extracolsep\fill}cccc} 
     \hline
 Nodes & $\epsilon$ & $\alpha$ & $\beta$ & Maximum Error \\
     \hline
9	 & 1.4440&	0.7404	&	0.0406 &	$5.79e-02$\\
16	 & 1.0864&	0.5993	&	0.0138 &	$1.59e-02$\\
25	 &1.1177&	0.6402  &	0.0239 &	$5.57e-03$\\
36	 &1.2702&	0.8400  &	0.2592 &	$3.05e-03$\\
49	 &1.3423&	0.7270  &	0.1930 &	$1.80e-03$\\
64	 &1.2631&	0.8148  &	0.2862 &	$7.76e-04$\\
81	 &1.0982&	0.9170  &	0.0287 &	$5.55e-04$\\
100	 &1.3937&	0.6212  &	0.1603 &	$5.41e-04$ \\
144	 &1.2657&	0.6009  &	0.2066 &	$2.05e-04$\\
196	 &1.1486&	0.7630  &	0.0256 &	$9.96e-05$\\
225	 &1.2052&	0.5768  &	0.0770 &	$8.00e-05$\\
400	 &1.1146&	0.7912  &	0.0175 &	$3.04e-05$\\
625	 &1.3459&	0.8404  &	0.2024 &	$1.22e-05$\\
900	 &1.2890&	0.5881  &	0.1275 &	$5.12e-06$\\
1600&1.5693&  0.8727  & 	0.2953 &	$1.56e-06$\\
2500&2.8189& 0.4630   &	0.1303 &	$6.47e-07$\\
     \hline
      \end{tabular*}
      \caption{The optimized values of parameters $\epsilon$, $\alpha$, and $\beta$ obtained during approximation of equation (\ref{possioneq}) at various degrees of freedom and corresponding maximum errors, in the numerical test \ref{s5.1}. The hybrid kernel performs better than both the Gaussian and cubic kernel.} 
      \label{NT1}
    \end{table*}
\subsection{Helmholtz Equation}
\label{s5.2}
\noindent Following Trefethen \citep{Trefethen2000}, who applied CHEB-PS method to numerically approximate 2D Helmholtz equation, and Fasshauer \citep{Fassbook2007}, who also used the same numerical example to test the implementation of RBFs in pseudospectral mode, we test the efficacy of the proposed HRBF-PS method for multivariate PDEs, using a similar test problem. 

Two-dimensional Helmholtz equation is given as 
\begin{eqnarray}
\label{Hequation}
\frac{\partial^2 u}{\partial x^2} +  \frac{\partial^2 u}{\partial z^2} + k^2 u = f(x,z), \qquad (x,z) \in (-1,1)^2,
\end{eqnarray}
with boundary conditions u=0, and
\begin{eqnarray}
f(x,z) = \exp\left( -10 \left[ \left(x-0\right)^2+\left(z-0\right)^2 \right]\right).
\end{eqnarray}
In order to solve this problem, we construct the Helmholtz operator on a tensor grid using \textit{Kronecker tensor-product}($\otimes$), as given by
\begin{eqnarray}
\label{HHoperator}
\mathcal{H} = \bm{D2}\otimes \bm{I}+ \bm{I}\otimes \bm{D2} + k^2\bm{I},
\end{eqnarray}
where $\bm{D2} = A_{xx}A_{xx}^{-1}$ is the second order differentiation matrix, and $\bm{I}$ is the identity matrix of size $(N+1)\times(N+1)$. It is to be noted that tensor-product grids are not necessary for RBF-PS algorithm. 
\subsubsection{Stability}
\noindent Table \ref{tab:NT2} enlists the condition number of various matrices arising during the numerical approximation of equation (\ref{Hequation}) using RBF-PS and HRBF-PS algorithms. $C_G$ and $C_H$ represent the condition numbers of the corresponding matrices, using Gaussian and hybrid kernel respectively. While implementing our modified algorithm, the parameters have been optimized with global PSO algorithm (see appendix A), using the cost vector obtained from LOOCV, as the objective function. The condition number of the interpolation matrix $A_{xx}$, second order differentiation matrices (D2), and the Helmholtz operators $\mathcal{H}$ are significantly reduced due to incorporation of the hybrid kernel in place of Gaussian kernel.  As observed from Table \ref{tab:NT2}, In case of Gaussian kernel, a sharp increase in the condition numbers is encountered after $24\times24$ collocation points in the domain, which gets smoother with hybrid kernel.

\begin{figure}[hbtp]
\centering
\includegraphics[scale=0.6]{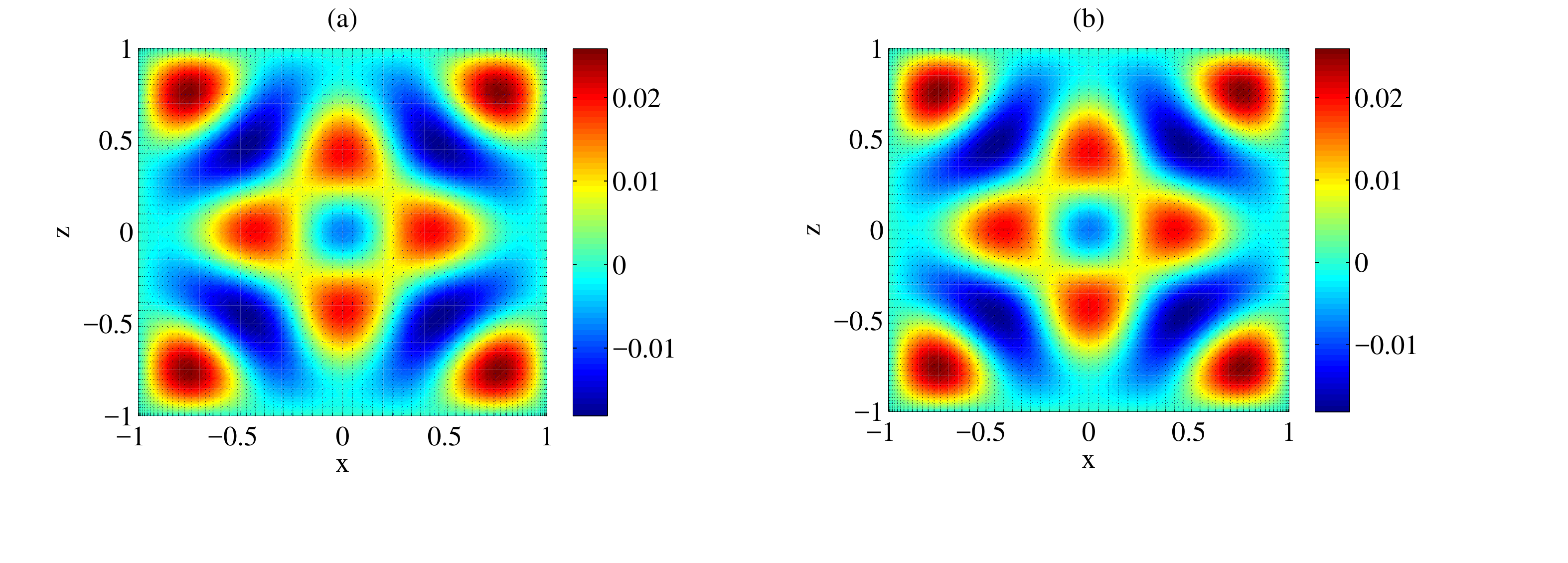}
\caption{The approximate solution of 2D Helmholtz equation with optimized parameters using (a) HRBF-PS and (b) CHEB-PS method, with 2500 nodes in the domain. Since the exact solution for this case was not known, the solution has been compared with Chebyshev pseudospectral method, which is known for its spectral accuracy. HRBF-PS with hybrid kernel shows excellent similarity with CHEB-PS.}
\label{fig:hybridvscheb}
\end{figure}

The stability of the operators is examined by plotting the corresponding eigenvalue spactra. It has been found that RBF based differentiation matrices have significant positive real part in the eigenvalue spectrum causing instability in the algorithm \citep{Palatte2006, Sarra2008}. Sarra \cite{sarra2011} proposed some numerical treatments to bound the condition number of the system matrix, to get stable eigenvalues of the RBF generated operators, which however, works with relatively large values of shape parameters. Here, we analyze the eigenvalues of the coefficient matrix constructed using both hybrid as well as Gaussian kernel, and compare them to those computed using CHEB-PS method. Eigenvalues for the low the degrees of freedom,\textit{ viz.}, upto $24\times 24$ collocation points in unit 2D domain, are stable for all the three methods. For higher degrees of freedom, however, RBF-PS algorithm becomes unstable due to emergence of some positive real eigenvalues in the spectra. This instability in RBF-PS algorithm increases with increasing degrees of freedom. The eigenvalue spectra of the HRBF-PS algorithm was found to be stable even at large degrees of freedom, which thus, mimics the same eigenvalue spectra, obtained using CHEB-PS approach. Since the exact solution of this problem is unknown, we compare the approximate solution obtained using HRBF-PS to spectrally accurate CHEB-PS method. As shown in Figure \ref{fig:hybridvscheb}, approximated solutions of this problem using HRBF-PS exhibits excellent agreement with that obtained by CHEB-PS approach.
\begin{table*}[!htbp]
\footnotesize
  \centering
  \begin{tabular*}{\textwidth}{l@{\extracolsep\fill}cccccc} 
    \hline 
    N & $C_G(A_{xx})$ & $C_H(A_{xx})$ & $C_G(D2)$ & $C_H(D2)$ & $C_G(\mathcal{H})$ & $C_H(\mathcal{H})$ \\
    \hline
   $ 5\times 5$   &$1.4335e+12$ & 26.1594 &$2.7930e+08$ & 342.4597 & 141.6702  & 196.3057 \\
   $ 9\times 9$   &$3.1443e+13$ & 158.3545 &$2.1967e+09$ & 757.6313 & $1.5601e+03$ & $1.7322e+03$ \\
   $13\times13$   & $1.9321e+13$ & 454.5354 &$2.8512e+09$ & $3.6177e+03$ &$7.6070e+03$ & $8.5408e+03$ \\
   $17\times17$   &$2.2286e+13$ & 990.7539  &$5.8481e+09$ & $1.1735e+04$& $2.4539e+04$ & $2.7889e+04$  \\
   $20\times20$   &$1.9861e+13$ & $1.5933e+03$ &$1.3034e+10$ & $2.4037e+04$ & $4.9853e+04$ & $5.7362e+04$ \\
   $24\times24$   &$1.5782e+13$ & $2.7210e+03$ &$1.3826e+10$ & $5.3886e+04$ &$ 1.1198e+05$ & $1.2915e+05$\\
   $36\times36$   &$5.1540e+17$ & $9.0045e+03$ &$3.5603e+13$& $3.2705e+05$ & $4.4562e+08$ & $7.9035e+05$\\
   $50\times50$   &$3.9606e+17$ & $2.3860e+04$ &$1.1026e+14$& $1.4176e+06$ & $1.1887e+09$ & $3.4435e+06$\\
   $64\times64$   &$3.0975e+18$ & $4.9726e+04$ &$8.8231e+14$& $4.2775e+06$ & $1.5040e+10$ & $1.0422e+07$\\
   $90\times90$   &$7.8193e+17$ &$1.0027e+05$ &$5.9061e+15$ &$1.2280e+07$ & $2.9728e+12$ & $3.0011e+07$\\
   $100\times100$ &$1.5220e+18$ &$1.8799e+05$ &$2.4384e+16$&$3.1582e+07$&$1.0992e+12$&$7.7335e+07$ \\  
\hline
\end{tabular*}
\caption{Condition number variation of the various matrices, at different degrees of freedom, in discretization of 2D Helmholtz equation via GRBF-PS and HRBF-PS}
\label{tab:NT2}
\end{table*}
\begin{figure}[hbtp]
\centering
\includegraphics[scale=0.35]{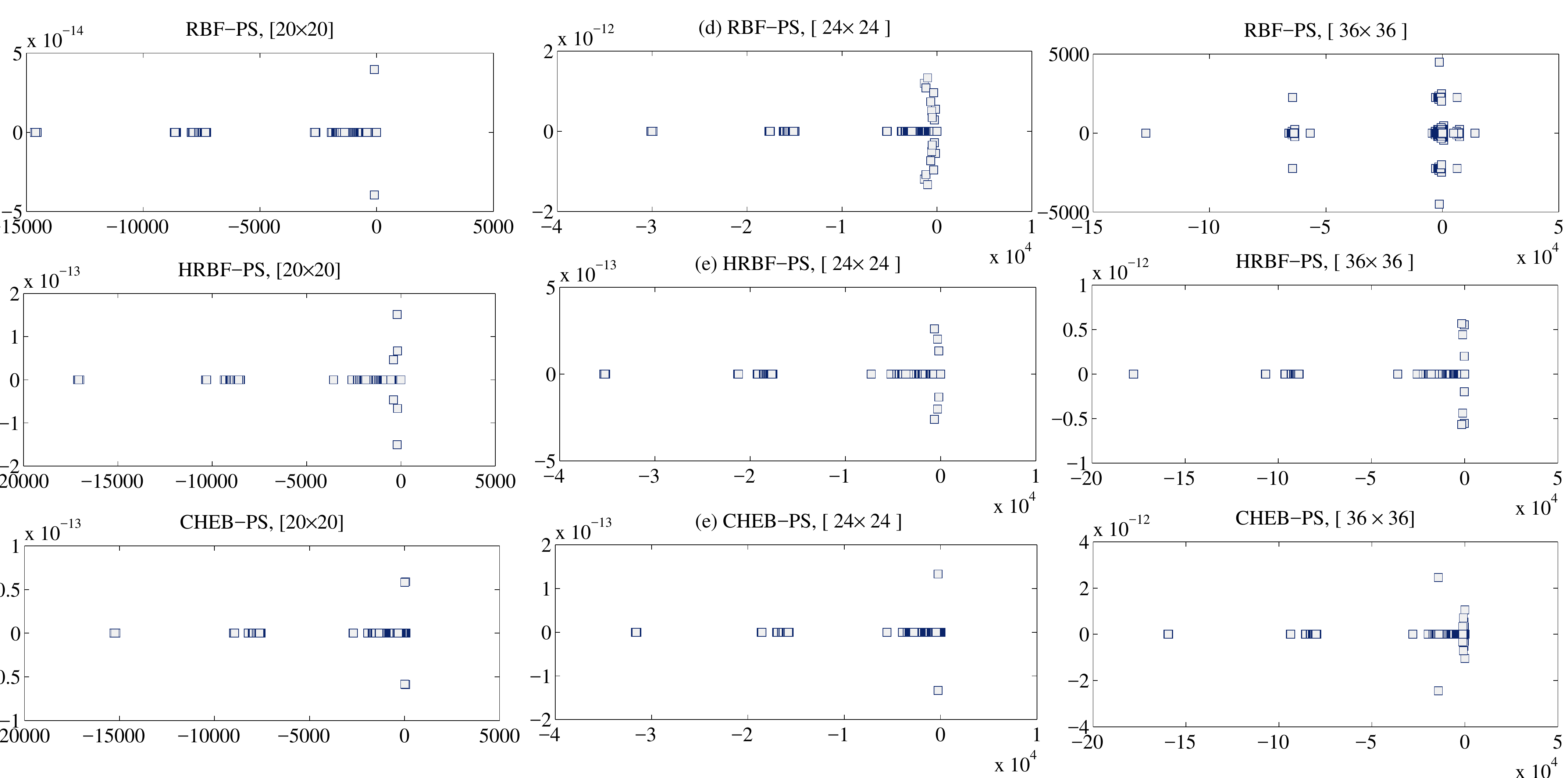}
\includegraphics[scale=0.35]{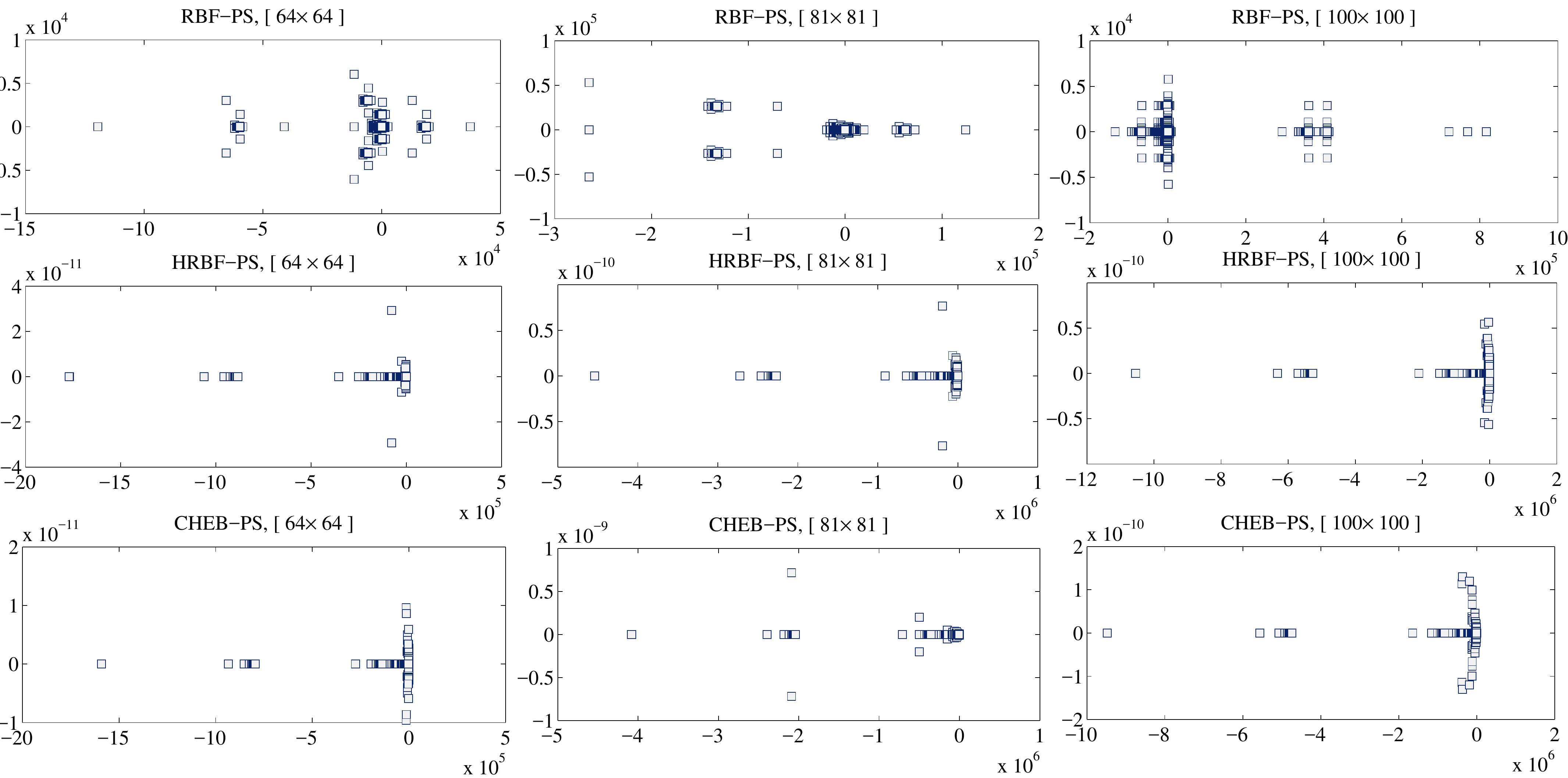}
\caption{Eigenvalue spectra of the coefficient matrices at different degrees of freedom in discretization of 2D Helmholtz equation via pseudospectral method using the Gaussian radial basis functions (GRBF-PS), hybrid radial basis functions (HRBF-PS) and Chebyshev polynomials (CHEB-PS) method. The real parts of the eigenvalues are on the horizontal axes, and the imaginary parts are on the vertical axes.}
\label{fig:eigen}
\end{figure}
\subsubsection{Computational Cost}
\begin{figure}[hbtp]
\centering
\includegraphics[scale=0.6]{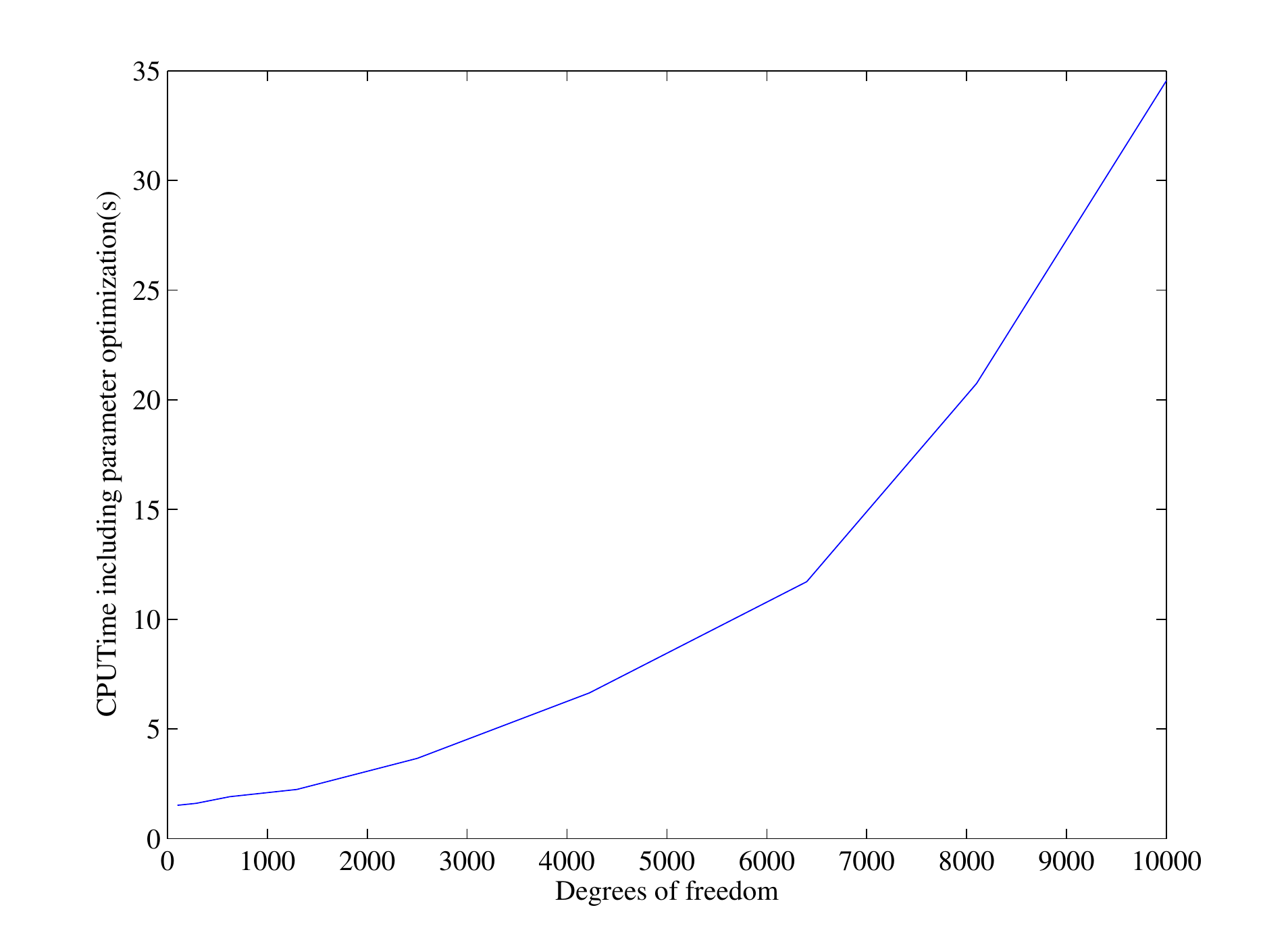}
\caption{Elapsed CPU time for solving 2D Helmholtz problem with RBF-PS scheme. This CPU time includes the time taken to find the optimized values of the parameters by PSO.}
\label{fig:CPU}
\end{figure}

\noindent In order to discuss the computational cost of the proposed approach, we measured the time taken by RBF-PS algorithm for various degrees of freedom. Table (\ref{tab:cost}) enlists the optimized values of parameters and corresponding elapsed CPU times, for solving the problem given by equation (\ref{Hequation}). This elapsed CPU time includes the time taken in optimization process. As shown in Figure (\ref{fig:CPU}), the cost of the present algorithm is roughly $O(N^3)$.
  \begin{table*}[htbp]
     \centering \footnotesize
     \begin{tabular*}{\textwidth}{l@{\extracolsep\fill}cccc} 
     \hline
 Nodes & $\epsilon$ & $\alpha$ & $\beta$ & CPU Time(s) \\
     \hline
100	 & 1.00&	$6.20e-01$	&	$2.06e-09$ &	1.52 \\
289	 & 1.21&	$7.69e-01$	&	$7.16e-08$ &	1.61 \\
625	 &1.37&	$8.07e-01$  &	$1.01e-06$ &   1.91 \\
1296 &1.00&  	$7.81e-01$  &	$5.33e-06$ &	2.24 \\
2500 &1.13&	$8.22e-01$  &	$1.82e-05$ &	3.66 \\
4225 &1.13&	$7.83e-01$  &	$5.26e-05$ &	6.64 \\
6400 &1.89&   $7.32e-01$  &	$1.23e-04$ &	11.71\\
8100&1.34&	    $8.78e-01$  &	$1.22e-04$ &	20.74\\
10000&1.29&  	$7.46e-01$  &	$2.00e-04$ &	34.55 \\
     \hline
      \end{tabular*}
      \caption{The optimized values of parameters $\epsilon$, $\alpha$, and $\beta$ obtained during approximation of equation (\ref{Hequation}) at various degrees of freedom and corresponding `Elapsed CPU time', in the numerical test \ref{s5.1}.} 
      \label{tab:cost}
    \end{table*}   
\begin{figure}[hbtp]
\centering
\includegraphics[scale=0.6]{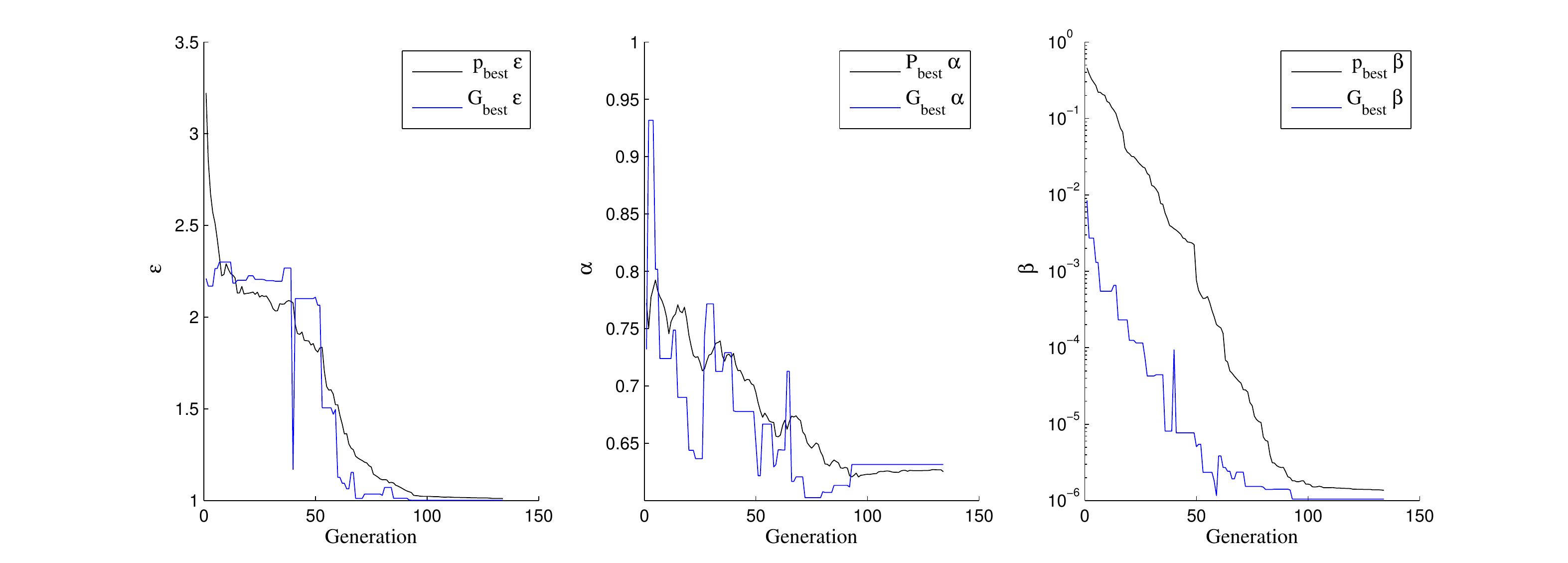}
\caption{The particle swarm optimization of the parameters for N=2500. The variation of (a) $\epsilon$ (b) $\alpha$ and (c) $\beta$ over 40 generations. }
\label{fig:acoustic2d}
\end{figure}
\subsubsection{Error variation with the shape parameter}
\begin{figure}[hbtp]
\centering
\includegraphics[scale=0.6]{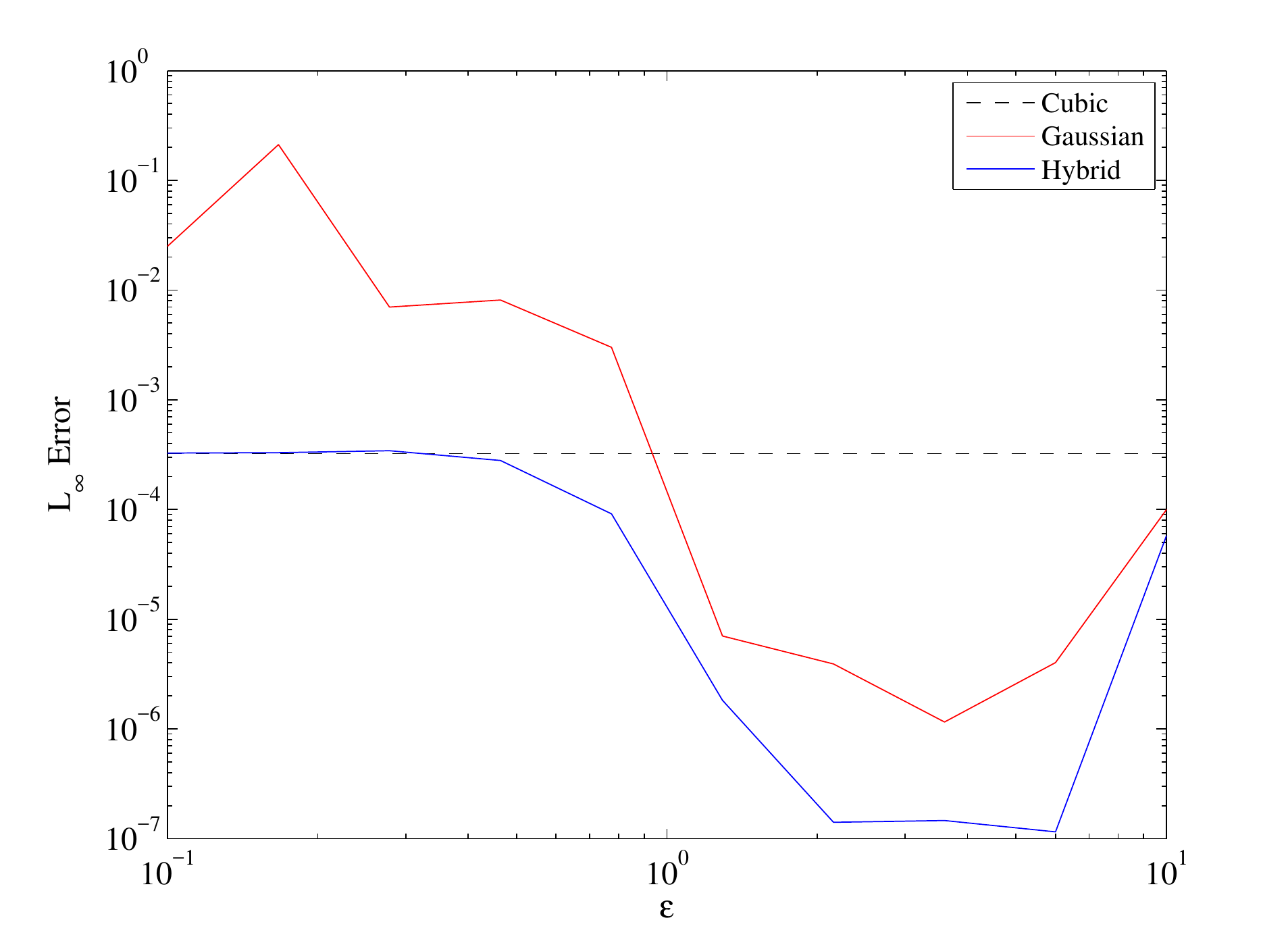}
\caption{Error variation with various values of the shape parameter $\epsilon$. This shows that the hybrid kernel ($\alpha =0.9, \beta =0.00001$) is more accurate than the Gaussian and cubic kernels.}
\label{fig:epprofile}
\end{figure}

\noindent Let us assume, a different case with equation (\ref{Hequation}), with the exact solution given by,
\begin{eqnarray}
\label{eq:HHexact}
u(x,z) = \frac{1}{1+x^2+z^2}.
\end{eqnarray}
The source term is given by, 
\begin{eqnarray}
f(x,z) = 8(x^2+z^2)u(x,z)^3-4u(x,z)^2+k^2u(x,z).
\end{eqnarray}
The Dirichlet boundary conditions are same as the exact solution. In order to compare the error variation with shape parameter of the Gaussian kernel, we fix the values of the weight coefficients as: $\alpha =0.9$ and $\beta =0.00001$, i.e., a very small doping of cubic kernel with the Gaussian. We compute the error using the approximated solution by RBF-PS method and the known exact solution given by equation (\ref{eq:HHexact}). As shown in Figure (\ref{fig:epprofile}), the hybrid kernel performs better than pure Gaussian and pure cubic kernel, however, in the case of very small shape parameter, the performance of the hybrid kernel, converges to the cubic. Figure (\ref{fig:epprofile}) also shows that with hybrid kernel and optimal value kernel parameters, RBF-PS method is more accurate than those with either of the two kernels. 
\subsection{Time-dependent PDE}
\noindent We applied HRBF-PS approach to solve a time dependent PDE by adopting a numerical example of 1D transport equation, from Fasshauer \citep{FASS20077}.

One-dimensional transport problem is given by
\begin{eqnarray}
\label{Tequation}
\frac{\partial u(x,t)}{\partial t} + c \frac{\partial u(x,t)}{\partial x}  = 0, \qquad x>-1,  t>0,
\end{eqnarray}
with the boundary conditions,
\begin{eqnarray}
\label{TBC1}
 u(-1,t)= 0, \nonumber
\end{eqnarray}
\begin{eqnarray}
\label{TBC2}
 u(x,0)= f(x). \nonumber
\end{eqnarray}
The analytical solution of this problem is expressed as 
\begin{eqnarray}
\label{TES}
 u(x,t) = f(x-ct). \nonumber
\end{eqnarray}

In order to compute the differentiation matrices using RBFs, Fasshauer \citep{Fassbook2007} suggested application of the Contour-Pade algorithm, which was initially proposed by Fornberg and Wright \citep{Forn2004853}. The Contour-Pad\'e algorithm allows the algorithm to evaluate radial basis function interpolants, with stability, for very small choices of the shape parameter, i.e., ``flat limit" ($\epsilon \rightarrow 0)$, which however, works only for small degrees of freedom. The results of Gaussian RBF and Chebyshev pseudospectral method has been taken from Fasshauer's \citep{Fassbook2007}.

\begin{figure}[hbtp]
\centering
\includegraphics[scale=0.4]{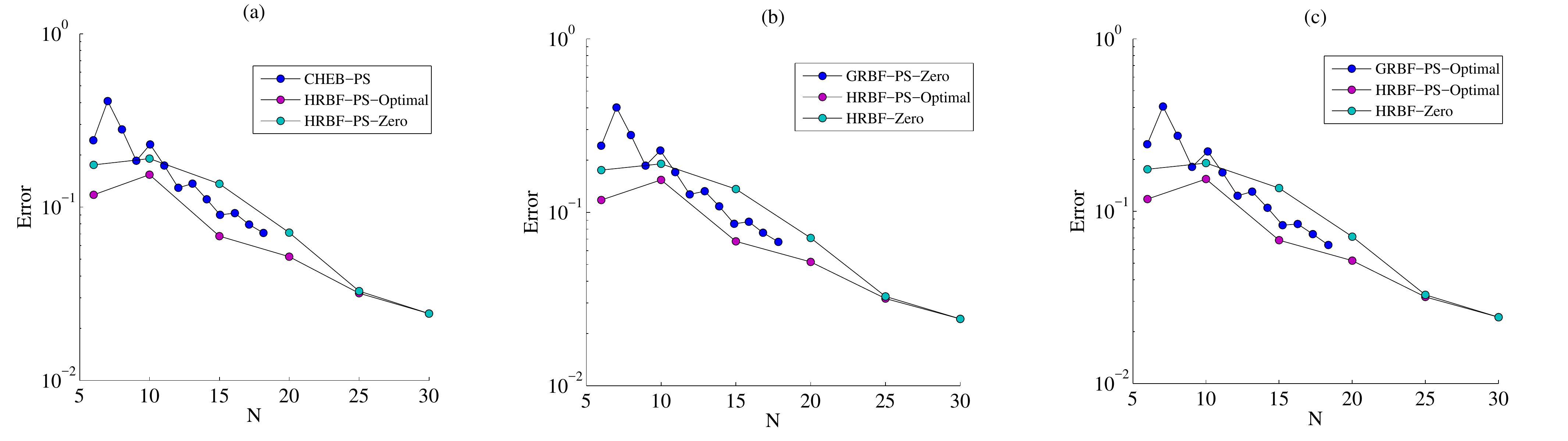}
\caption{Numerical approximation of transport problem: the convergence patterns obtained using the hybrid kernel with  $\epsilon =0$ and \textit{optimal $\epsilon$}, have been compared to that obtained with Chebyshev pseudospectral method (a), the Gaussian RBF with $\epsilon=0$ (b), and the Gaussian RBF with \textit{optimal} $\epsilon$ (c). The data plotted with blue lines have been taken from Fasshauer's book \citep{Fassbook2007}, which suggest that the errors in this numerical tests are mostly due to the time-stepping method.}
\label{fig:transport}
\end{figure}

In the diagrams depicted in Figure \ref{fig:transport}, the maximum errors at time $t=1$ with time step $\Delta t =0.001s$ (using implicit Euler method for time stepping) have been shown. According to Fasshauer \citep{Fass06,Fassbook2007}, Contour-Pad\'e allows a limited spatial discretization (upto $N=18$ only), in this case. With the hybrid kernel, however, we could use relatively large number of nodes. The convergence patterns obtained using the HRBF-PS algorithm with zero and optimized shape parameters, have been compared to those obtained using CHEB-PS (Figure \ref{fig:transport}a), and RBF-PS algorithms (Figure \ref{fig:transport}c). It can be seen that the proposed hybrid kernel not only improves the convergence but also numerically solves this time-dependent PDE beyond the limitation imposed by Contour-Pad\'e approach. As shown in \citep{Fassbook2007}, in this context, PS and RBF-PS show ``virtually identical" error pattern, therefore, the proposed hybrid kernel exhibit similar improvement in all the three plots.
\subsection{Laplace equation with non-trivial boundary conditions}
\label{sec:nt4}
\begin{figure}[hbtp]
\centering
\includegraphics[scale=0.45]{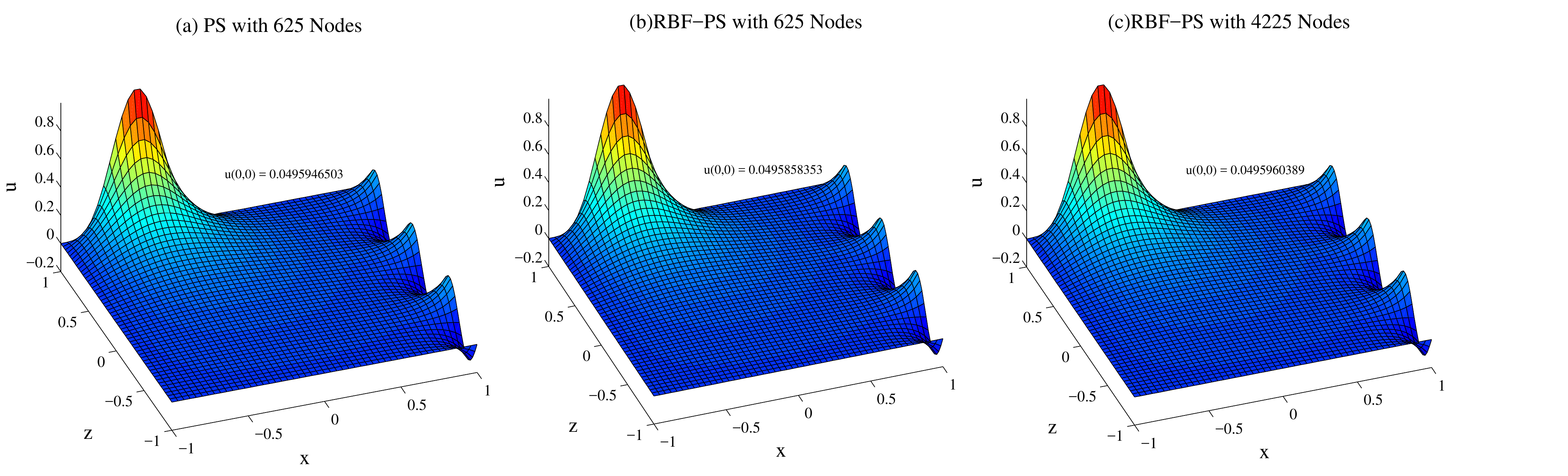}
\caption{The solution of the problem, considered in the section [\ref{sec:nt4}] (a) using PS method with 625 nodes, (b) using RBF-PS method with 625 nodes, and (c) using RBF-PS method with 4225 Nodes. The hybrid kernel has been used in the RBF-PS formulation. The optimal values of the parameters $[\epsilon, \alpha, \beta]$ were decided by the particle swarm optimization, which are $[1.2652, 0.81219, 4.7835e-05]$ and $[1.0476, 0.76712, 4.6704e-05]$ for cases (b) and (c), respectively.}
\label{fig:lap2d}
\end{figure}

\noindent In order to test the efficacy of the proposed RBF-PS method, we consider another example from Trefethen\citep{Trefethen2000}, as given by
\begin{eqnarray}
\frac{\partial^2 u}{\partial x^2} +  \frac{\partial^2 u}{\partial z^2} = 0, \qquad x,z \in (-1,1)^2,
\label{eq:lap2d}
\end{eqnarray}
\noindent with piece-wise boundary conditions
\begin{eqnarray}
   u(x,z)= 
\begin{cases}
    \sin^4(\pi x), & \text{if } z=1, -1< x < 0,\\
    \frac{1}{5}\sin(3\pi x), & \text{if } x=1, \\
    0,              & \text{otherwise}.   
\end{cases}
\label{eq:lap2dbc}
\end{eqnarray}
\noindent Figure \ref{fig:lap2d} shows the solution of above problem with the presented RBF-PS method. We optimize the kernel parameters considering that exact solution of the problem was not known and compare the solution with that obtained using spectrally accurate PS method. The solutions obtained by using RBF-PS method exhibit excellent similarity with that obtained by using PS method, which gets closer to the accuracy of PS method at higher degrees of freedom. This numerical test shows that the presented hybrid kernel can handle the piece-wise non-zero boundary conditions. Also, the presented approach not only keeps the algorithm stable at relatively higher degrees of freedom but also provides excellent accuracy.
\section{Conclusion}
\noindent We proposed a novel implementation of hybrid Gaussian-cubic kernels in radial basis-pseudospectral approach for numerical approximation of PDEs.  Such hybrid kernels make the algorithm well-posed and enable it to perform computation with relatively larger degrees of freedom as well as with very small shape parameters. Based on the numerical tests performed in this study, we draw the following conclusions.
\begin{enumerate}
\item Application of pure Gaussian kernel in an RBF-PS algorithm leads to ill-conditioning of the resulting linear system, on the other hand, using purely spline kernels will create the risk of singularity for certain node arrangements. The presented hybrid kernel, is therefore, a reasonable choice to be used in global RBF-PS algorithms. 

\item The accuracy of the hybrid kernels is found to be better than pure Gaussian as well as polyharmonic spline (cubic) kernel. Also, with hybrid kernel, the global RBF-PS algorithm does not diverge at higher degrees of freedom, ensuring the stability of the proposed approach.

\item The eigenvalue spectra of the system matrix were found to be stable at relatively larger degrees of freedom, which ensures the stability of the linear system. 

\item Finding the ``optimal value" of the shape parameter for kernel-based meshless algorithms is conventional, however, the two extra parameters introduced due to hybridization of two kernels is more likely to increase the computational cost. The cost of the present approach is found to vary as $O(N^3)$, which is similar to RBF-QR approach. However, the cost of RBF-QR increases drastically with increase in shape parameter \citep{Forn2004853}.

\item Using particle swarm optimization for finding the parameter(s) of a RBF in kernel-based meshless computing is a novel approach, recently proposed in \citep{1512.07584}.
Here, we provide a detailed discussion about the application of PSO  for finding the shape parameter and weight coefficients of RBFs for kernel based meshless computing (see Appendix A), for numerical approximation of PDEs. This algorithm can be easily simplified for shape parameter tuning in other meshless algorithms.

\item We have tested the improvements in RBF-PS method with the presented approach, and shown them using similar numerical tests, used in the initial development of RBF-PS method. Further studies can explore the application of the presented approach for more complex problems, and with different node discretization approaches.
\end{enumerate}
\newpage
\section*{Appendix A: Particle Swarm Optimization}
\noindent PSO is a powerful optimization technique which uses the intelligence of swarms to solve problems. It was developed in 1995 by Kennedy and Eberhart \citep{Eberhart1995}. In PSO, a certain number of assumed solutions, termed as \textit{particles} are initialized and further directed towards the best solution over a defined number of generations. Each particle is a point in the search space which ``flys" using its own experience as well as the experience of other \textit{particles}. The position of each \textit{particle} in the solution space is tracked by the algorithm. The \textit{particle} associated with the best fitness so far is termed as its personal best value \textit{pbest}. The best value obtained so far by any particle at a generation is termed as global best value \textit{gbest}. The idea behind PSO algorithm is to accelerate each particle towards its \textit{pbest} and \textit{gbest}. In the discussion below, we will focus on the PSO, specifically applied to find the shape parameter and weight coefficients, for the hybrid kernel proposed in this paper. However, for general purpose understanding of PSO, we recommend \citep{Marini2015,Singh2015}.

\begin{figure}[hbtp]
\centering
\includegraphics[scale=0.55]{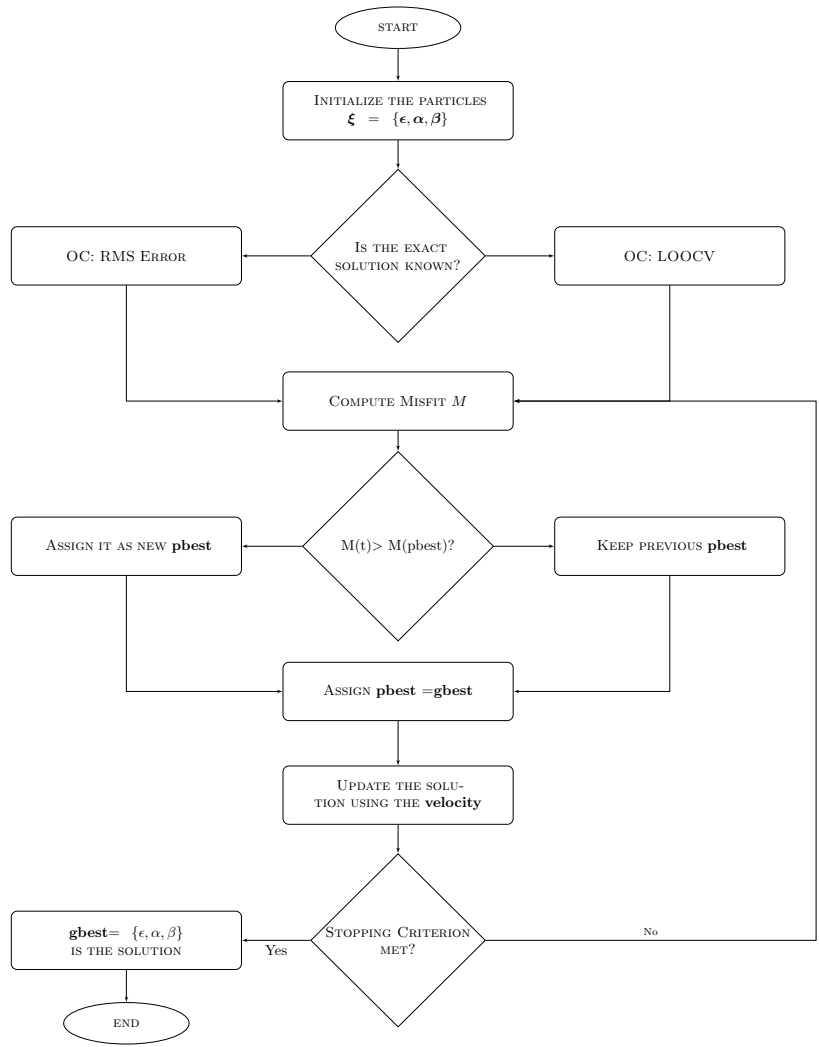}
\caption{A flowchart showing particle swarm algorithm in the context of finding $\epsilon$, $\alpha$, and $\beta$. The optimization criterion (OC) is decided using RMS error or the cost function provided by leave-one-out-crossvalidation (LOOCV), depending upon the availability of the exact solution.}
\label{fig:psoflowchart}
\end{figure}

Let us assume a particle $\bm{\xi} = \{\bm{\epsilon}, \bm{\alpha}, \bm{\beta}\}$. This particle $\bm{\xi}$ contains certain number of values of each parameters within the user specified search ranges. The position of the $i^{th}$ particle can be written as,

\begin{eqnarray}
\xi_i = [\epsilon_i, \alpha_i, \beta_i]. 
\end{eqnarray}
Until the stopping criterion is met, the position of this particle in solution space is continuously updated according to the following equation,
\begin{eqnarray}
\xi_i(t+1) = \xi_i(t) + v_i(t+1) 
\end{eqnarray}
where t indicates a typical iteration of the algorithm and $v_i$ represents the vector collecting the velocity-components of the $i^th$ particle. This velocity vector decides the movement strategy of the particle in the search space. At each iteration, the velocity component is updated according to the equation given below,
\begin{eqnarray}
\label{eq:velocity}
v_i(t+1) = v_i(t) + c_1\left(p_i-\xi_i(t)\right)R_1 + c_2\left(g-\xi_i(t)\right)R_2. 
\end{eqnarray}
There are three components of the velocity vector, which are represented by three terms on the right hand side of equation (\ref{eq:velocity}). The first term is defined as \textit{inertia}, which prevents substantial changes in the particle direction. The second term is defined as the \textit{cognitive-component}, which controls the particles' tendency to return to their own previously found best solutions. Finally, the third term is defined as the \textit{social-component}, which quantifies the best particle relative to its neighbors. The best solution obtained so far, in an iteration, by a specific individual is termed as ``personal best" ($p_i$). The ``global best" (g), on the other hand, represents the overall best solution obtained by that particular swarm. The ``cognitive coefficient'' $(c_1)$ and ``social coefficient'' ($c_2$) are real-valued constants, which moderates, in each iteration, the steps taken by a particle in the direction of its ``personal best" and ``global best", respectively. Perez and Behdinan \cite{Perez2007} proposed a constraint as $0< c_1+c_2 < 4$, to ensure the stability in the global PSO algorithm. We have kept the values as $c1=1.2$ and $c_2=1.7$, however, it was observed that in this context, the optimization works fine as long as $0< c_1+c_2 < 4$ is satisfied. Figure (\ref{fig:psoflowchart}) shows the PSO algorithm, applied in the context of this paper. \\

\newpage
\bibliographystyle{elsarticle-num}
\bibliography{sample}
\end{document}